\documentclass[11pt]{article}
\usepackage[small,compact]{titlesec}
\usepackage{amssymb, amsmath, amsthm, ulem}

\begin{document}

\thispagestyle{empty}
\begin{center}
\Large
Authors' Note
\end{center}

This is a corrected version of the article

M. E. H. Ismail and M. E. Muldoon, {\it Inequalities and monotonicity properties for gamma and $q$-gamma functions}, pp. 309-323 in R. V. M. Zahar, ed., Approximation and Computation: A Festschrift in Honor of Walter Gautschi, ISNM, vol. 119, Birkh\"auser, Boston-Basel-Berlin, 1994. 

Most of the errors in the original paper had to do with saying that certain functions related to 
the $q$-gamma function were {\bf not} completely monotonic.  We discovered these errors through reading the paper {\it Some completely monotonic functions involving the $q$-gamma function}, by Peng Gao, \\{\tt http://arxiv.org/abs/1011.3303}.

We also take the opportunity to correct some errors in other places including the statement and proof of Theorem 3.4.
\newpage
 \setcounter{page}{309}
 \noindent {\bf Corrected July 31, 2011}
\begin{center}
\LARGE{\sc  Inequalities and monotonicity properties \\for gamma
and $q$--gamma functions } 

\normalsize
Mourad E. H.  Ismail\footnote{Department of Mathematics,
University of South Florida, Tampa, FL 33620--5700, U. S. A.}
\mbox{  }
Martin E. Muldoon\footnote{
Department of Mathematics and Statistics, York
University, North York, Ont. M3J 1P3, Canada}
\\
\bigskip
{\sl To Walter Gautschi on his
65th birthday}
\end{center}
\begin{quotation}
{\bf Abstract}
We prove some new results and unify the proofs of old ones
involving complete monotonicity of expressions involving gamma
and $q$--gamma functions, $0 < q < 1$. Each of these results
implies the
infinite divisibility of a related probability measure.  In a few
cases, we are able to get simple monotonicity without having
complete monotonicity.  All of the results lead to inequalities for
these functions.  Many of these were motivated by the bounds in a
1959 paper by Walter Gautschi.   We show that some of the bounds
can be extended to complex arguments.
\end{quotation}


\section{INTRODUCTION AND PRELIMINARIES}
\setcounter{equation}{0}
Among Walter Gautschi's many contributions to mathematics are some
interesting inequalities for the
gamma function.  For example, he shows (\cite{gautschi74a},
\cite{gautschi74b}) that 
if $x_k >0,k=1,\dots n,\; x_1x_2\cdots x_n =1$, then the
inequality
\begin{equation} \sum_{k=1}^n \frac{1}{\Gamma(x_k)} \le n
\end{equation}
is true for $n = 1,2\dots, 8$ but not for $n \ge 9$. 
 Here we will be more concerned with an earlier result of
Gautschi's \cite{Ga}, the two-sided inequality
\begin{equation}
n^{1-s} < \frac{\Gamma(n+1)}{\Gamma(n+s)} < \exp[(1-s)\,
\psi(n + 1)],\quad
0 < s < 1,\;n=1,2,\dots 
\end{equation}
which still inspires extensions. For example, D. Kershaw
\cite{kershaw} proved
\begin{eqnarray}
\exp[({1-s}) \psi(x + s^{1/2})] < \frac{\Gamma(x+1)}{\Gamma(x+s)}
& < & \exp[(1-s)\,
\psi(x+(s+1)/2)], \nonumber \\
& & 0 < s < 1, \; x > 0, 
\end{eqnarray}
and
\begin{equation}
\left[x + \frac{s}{2} \right]^{1-s} <
\frac{\Gamma(x+1)}{\Gamma(x+s)} < 
\left[x -\frac12 + \left(s + \frac14 \right)^\frac{1}{2}
\right]^{1-s}, \;\;
0 < s < 1,\; x> 0. 
\end{equation}
\newpage
\noindent In all of these inequalities, $\psi(x)$ denotes the logarithmic
derivative $\Gamma'(x)/\Gamma(x)$. 

 Many inequalities for special functions follow from monotonicity
properties. Often such inequalities are
special cases of the complete monotonicity of related special
functions. For example, an inequality of the form $f(x)\ge g(x)$,
$x\in [a,\infty)$ with equality if and only if $x=a$, may be a
disguised  form of the complete monotonicity of
$g(\varphi(x))/f(\varphi(x))$ where $\varphi$
 is a nondecreasing function on $(a,\infty)$ and
$g(\varphi(a))/f(\varphi(a))=1$.
Thus, for example, the left--hand inequality in (1.2) and
the right--hand one in (1.3) follow from the
facts that $x^s \Gamma(x+s)/\Gamma(x+1)$ and 
$\exp[(s-1)\,
\psi(x+(s+1)/2)] \Gamma(x+s)/\Gamma(x+1)$ are, respectively,
decreasing and increasing functions of $x$ on $(0,\infty)$. 
Bustoz and Ismail \cite{Bu:Is} proved that some of the above
inequalities for the gamma function follow from the complete
monotonicity of certain functions involving the ratio
$\Gamma(x+1)/\Gamma(x+s)$. 

Recall that a  function $f$ is completely monotonic on an
interval $I$ if  $$(-1)^n f^{(n)} (x) \ge 
   0$$ for  $n = 1,2,\dots$ on $I$. We collect some known facts,
all
either easily proved or contained in \cite{Wi}  or \cite{Fe},  in
the following theorem.

{\bf Theorem 1.1. } {\sl
(i)  A necessary and sufficient condition that $f(x)$ should be
completely monotonic on $(0,\infty)$ is that 
$$f(x) = \int_0^\infty e^{-xt} d\alpha(t),$$
where $\alpha(t)$ is nondecreasing and the integral converges for
$0 < x < \infty$.

(ii) $e^{-h(x)}$ is completely monotonic on $I$ if
$h'(x)$ is  completely monotonic on I. 

(iii) A probability  distribution    supported on a subset of
$[0,\infty)$ is infinitely divisible if and only  if its Laplace
transform (moment generating function) is of the form
   $e^{-h(x)}$ with $h(0^+) = 0$ and $h'(s)$ is completely
monotonic on 
   $(0,\infty)$.
}

For brevity, we shall use {\sl completely monotonic} to mean {\sl
completely   monotonic on $(0,\infty)$}.

There is already an extensive and rich literature on inequalities
for gamma
functions; for references see \cite{Ma:Ol}, \cite{Mi}. 
One of the objects of the present work is to show that many of
these can be extended, using essentially the same methods of proof,
to the $q$--gamma function  defined (see, e.g., \cite{Ga:Ra}) by 
\begin{equation}
\Gamma_q(x) := (1-q)^{1-x}\prod_{n = 0}^\infty \frac{1-q^{n+1}}{1-
q^{n+x}},\quad
 0 < q < 1. 
\end{equation}
Although the right-hand side of (1.5) is meaningful when
$|q| < 1$, our results will require $q^x \in (0, 1)$ for all
positive $x$.
This forces $q\in (0, 1)$. As expected
\begin{equation}
\Gamma_q(x)\to \Gamma(x) \quad as \; q\to 1^-.
\end{equation}
\newpage
The most elegant proof of this, due to R. William Gosper, is
 in Appendix A of Andrews's excellent monograph \cite{An}; see also
\cite[p. 17]{Ga:Ra}.  For a rigorous justification, see
\cite{koorn}.   It is worth noting that
\begin{equation}
\Gamma_q(x) \approx (1 - q)^{1-x} \prod_{0}^{\infty}(1 -
q^{n+1}),\quad  {\rm as}\; x \to
\infty\;{\rm  if}\; |q| < 1.
\end{equation}
It seems that most of our results have analogues also for the
$q$--gamma function with $q > 1$, in which case the definition
(1.5) must be changed.  We do not pursue this question here.
For the gamma function, we have the (Mittag--Leffler) sum
representation \cite{Er:Ma1}
\begin{equation}
\psi(z) = \frac{\Gamma'(z)}{\Gamma(z)} = -\gamma +
 \sum_{n=0}^\infty \left(\frac{1}{n+1} - \frac{1}{z+n}\right), 
\end{equation}
and the integral representation \cite[(1.7.14)]{Er:Ma1}
\begin{equation}
\psi(x) = -\gamma + \int_{0}^\infty \; \frac{e^{-t}- e^{-tx}}{1 -
e^{-t}} \, dt,\quad  {\rm Re} \;x > 0. 
\end{equation}
Although (1.9) and (1.8) are equivalent, for
${\rm Re} \;x > 0$, it turns
out that (1.9) is more useful in proving the kind of 
inequalities in which we are interested.  This situation occurred
also in \cite{ronning}.

A corresponding sum representation for the case of the $q$--gamma
function, easily following from (1.5) is
\begin{equation}
\psi_q(x) := \Gamma'_q(x)/\Gamma_q(x) = - \log  (1-q) +  \log  q
\sum_{n=0}^\infty
q^{n+x}/(1 - q^{n+x}),\;0 < q <1.  
\end{equation}
Although this representation has been used directly in the proofs
of many results for the $q$--gamma function (in \cite{Bu:Is}, for
example) we shall find it more convenient to use the equivalent 
Stieltjes integral representation
\begin{equation}
\psi_q(x) = - \log  (1-q) -  \int_0^\infty
\frac{e^{-xt}}{1-
e^{-t}} d\gamma_q(t), \; 0 < q < 1,\; x >0, 
\end{equation}
where $ d\gamma_q(t)$ is a discrete measure with positive masses
$-\log q$ at the positive points $-k \log q,\;\;k=1,2,\dots.$
For completeness, and economy of later statements, we include the
value $q =1$ in the definition of $ \gamma_q(t)$:
\begin{equation}
 \gamma_q(t) = \left\{
\begin{array}{l} -\log q \sum_{k=1}^\infty \delta(t + k \log q),\;
0< q <1,\\
t,\;\;q=1.
\end{array} \right. 
\end{equation}
To get the representation (1.11), we expand the
denominator of the sum in
(1.10) by the binomial theorem and interchange the orders 
of summation to get
\begin{equation}
\psi_q(x) = - \log  (1-q) +  \log  q \sum_{k=1}^\infty q^{kx}/(1 -
q^{k}),  
\end{equation}
\newpage
\noindent which is equivalent to (1.11).  

Note that $\psi_q(x)$ can also be expressed as a $q$--integral
\cite[p. 19 ]{Ga:Ra},
\begin{equation}
\psi_q(x) = - \log  (1-q) +  \frac{\log  q}
{1-q} \int_0^1 \frac{t^{x -1}}{1 - t} d_q(t),\;0 < q< 1,\;x >0,
\end{equation}
just as, from  (1.9), $\psi(x)$  can be expressed as an
ordinary integral over $[0,1]$:
\begin{equation}
\psi(x) = -\gamma + \int_{0}^1\; \frac{1 - t^{x-1}}{1 - t}
\, dt,\quad  {\rm Re} \;x > 0. 
\end{equation}
We will need the following relations which follow easily from the
definition of $d\gamma_q(t) $:
\begin{equation}
\int_0^\infty e^{-xt} d\gamma_q(t) = \frac{-q^x \log q}{1 -
q^x},\;0 < q<
1,\;x >0,
\end{equation}
and
\begin{equation}
 \int_0^\infty \frac{e^{-xt}}{t} d\gamma_q(t) = 
\sum_{k=1}^\infty \frac{q^{kx}}{k} = - \log(1-q^x),\;0 < q< 1,\;x
>0.
\end{equation}
We will use the following Lemma in many of our proofs. We remark
that it includes results stated in different notations 
\cite[Lemma 3.1]{Bu:Is} and \cite[Lemma 4.1]{Is:Lo},
as well as individual steps proved by {\sl ad hoc} methods in these
and other papers. 

{\bf Lemma 1.2.} {\sl
Let $0 < \alpha < 1$. Then
\begin{equation}
\alpha e^{(\alpha -1)t} < \frac{{\rm sinh}\, \alpha t} {{\rm sinh}\, t}
< \alpha,\;\;\; t > 0.
\end{equation}
The inequalities become equalities when $\alpha = 1$ and they are
reversed when $\alpha > 1$.
}

The following lemma will also be useful.

{\bf Lemma 1.3.} {\sl (i) Let $f(x)$ be completely monotonic on
$(0,\infty)$
and
let $a > 0$.  Then $f(x) - f(x +a)$ is completely monotonic on
$(0,\infty)$. (ii) Let $f(x) \ge 0$ and let $f(x) - f(x +a)$ be
completely monotonic on $(0,\infty)$ for each $a$ in some 
right--hand neighbourhood of $0$. Then $f(x)$ is completely
monotonic on $(0,\infty)$. 
}

{\bf Proof:} 
(i) We have
$$f(x) = \int_0^\infty e^{-xt} d\alpha(t),$$ 
where $\alpha(t)$ is nondecreasing and the integral converges for
$0 < x < \infty$.
Hence 
$$(-1)^n D_x^n [f(x) -f(x +a)] = \int_0^\infty [e^{-xt} - e^{-
(x+a)t}] t^n d\alpha(t) \ge 0.$$ 
(ii) Under the hypotheses here, we find that 
$-f'(x) = \lim_{a \rightarrow 0^+} [f(x) - f(a+x)]/a$ is completely
monotonic on $(0,\infty)$. \hfill $\blacksquare$
\newpage
{\bf Remark 1.} A feature of the present work is that the
similarity between 
(1.9) and (1.11) makes it possible to unify the
proofs of some of our results for the gamma and $q$--gamma
functions.

{\bf Remark 2. } The integral representation in
Theorem 1.1 (i) provides a necessary as well as a sufficient
condition for the complete monotonicity of $f$. This enables us to
show that certain monotonic functions are {\sl not} completely
monotonic more easily than is done in \cite{alzer}, for example. As
in \cite{Is:Lo}, many of our results will assert the complete
monotonicity of a function for a certain range of values of a
parameter, the complete monotonicity of its derivative for another
range, and, {\bf in some cases}, the complete monotonicity of neither of these for an
intermediate range.

{\bf Remark 3.} In many discussions of completely monotonic
functions, the concept of {\sl strict complete monotonicity} is
used to indicate strict inequality in $(-1)^n f^{(n)} (x) \ge 
0$.  But if, as here, our interval is a half--line, we get such
strict inequality in all but trivial cases:  A result of
J. Dubourdieu \cite[p. 98]{dub} asserts that for a completely
monotonic function on $(a,\infty)$, we have $(-1)^n f^{(n)} (x) >
0$ for  $n = 1,2,\dots$, unless $f(x)$  is constant.

{\bf Remark 4.} In \S5, we extend some bounds for ratios of gamma
functions to complex values of the arguments.



\section{GAMMA AND $q$--GAMMA FUNCTIONS}
\setcounter{equation}{0}

The following result was proved in \cite{Is:Lo}: 

{\bf Theorem 2.1.} {\sl
Let $h_\alpha(x) = \log[x^\alpha\Gamma(x) (e/x)^x]$. Then $-
h_\alpha{'}(x)$ is completely monotonic on $(0,\infty)$ for 
$\alpha \le 1/2$, $h_\alpha{'}(x)$ is completely monotonic 
for $\alpha \ge 1$, and neither is completely
monotonic for $1/2 < \alpha < 1$.
}

The proof in \cite{Is:Lo} is based on the consequence 
$$ -h_\alpha'(x) = \int_0^\infty \left[\frac{1}{1-e^{-t}} -
\frac{1}{t} -\alpha \right] e^{-xt} dt, $$ of  (1.9)
and the fact that the quantity in the square brackets, which has
the same sign as
$(1 +\alpha t)e^{-t} -1 + (1 -\alpha t)$, is positive for $\alpha
\le 1/2$, negative for $\alpha \ge 1$ and undergoes a change of
sign for $1/2 < \alpha < 1$.
The next result can be considered a $q$-analogue of Theorem 2.1.

{\bf Theorem 2.2.} {\sl
Let $ 0< q < 1$ and let
$$h_\alpha(x) = \log\left[
(1-q)^x(1-q^x)^\alpha\Gamma_q(x) 
\exp\left({\sum_{k=1}^\infty q^{kx}/(k^2 \log q)}\right) \right].$$
Then $-h_\alpha{'}(x)$ is completely monotonic on $(0,\infty)$ for 
$\alpha \le 1/2$ {and} $h_\alpha{'}(x)$ is completely monotonic on
$(0,\infty)$ for $\alpha \ge 1$.
}
 \newpage
{\bf Proof:} 
It follows from (1.11) that
$$ - h_\alpha'(x) =  \int_0^\infty\left[\frac{1}{1-e^{-t}} -
\frac{1}{t} -\alpha \right] e^{-xt} d\gamma_q (t), $$
where $d \gamma_q(t)$ is defined by (1.12).  As in
the proof of \cite[Theorem 2.1]{Is:Lo}, the quantity in square
brackets is positive for $\alpha \le \frac12$ { and} negative for
$\alpha \ge 1$. 
$\blacksquare$

{\bf Remark.}   To see that Theorem 2.1  { includes} the limiting case of
Theorem 2.2 as $q
\rightarrow 1^-$, we will show that 
$$
\lim_{q \rightarrow 1^-}
\frac{(1-q)^x(1-q^x)^\alpha\Gamma_q(x) \exp [ F(q^x)/\log q]}
{(1-q)^\alpha
\exp[F(1)/\log q]} = x^\alpha \Gamma(x)(e/x)^x,
$$
where 
\begin{equation} F(x) = \sum_{n=1}^\infty \frac{x^n}{n^2} = -
\int_0^x
\frac{\log(1-t)}{t} dt.
\end{equation}
There is no difficulty in seeing that
$$
\lim_{q \rightarrow 1^-}
\frac{(1-q^x)^\alpha\Gamma_q(x) }{ (1-q)^\alpha} = x^\alpha
\Gamma(x),
$$
so it remains to show  that 
$$
\lim_{q \rightarrow 1^-}
\frac{(1-q)^x\exp[F(q^x)/\log q]} {\exp[ F(1)/\log q]} = (e/x)^x.
$$
Taking logarithms, this is equivalent to showing that 
$$
\lim_{q \rightarrow 1^-} \left[x \log (1-q) + \frac{F(q^x) -
F(1)}{\log q}\right]
=x - x \log x,
$$
and this follows from the identity
$$x \log (1-q) + \frac{F(q^x) - F}{\log q} \equiv - \int_0^x
\log
\frac{1-q^t}{1-q} dt,
$$
which is easily checked by differentiation with respect to $x$.

Applying Lemma 1.3 to the results of the last theorem, we get:

{\bf Theorem 2.3.} {\sl
Let $ 0< q < 1$, $a >0$ and let
$$H_\alpha(x) = \log\left[
\left( \frac{1-q^x}{1-q^{x+a}} \right)^\alpha
\frac{\Gamma_q(x)}{\Gamma_q(x+a)}
\exp\left(\sum_{k=1}^\infty \frac{q^{kx} - q^{k(x+a)}}{k^2 \log
q}\right) \right].$$
Then $-H_\alpha{'}(x)$ is completely monotonic on $(0,\infty)$ for 
$\alpha \le 1/2$ {and}  $H_\alpha{'}(x)$ is completely monotonic on
$(0,\infty)$ for $\alpha \ge 1$.
}
\newpage
In the case $q \rightarrow 1^-$,  the following
Corollary  follows by applying Lemma 1.3 to the result
of Theorem 2.1. 

{\bf Corollary 2.4.} {\sl
Let  $a >0$ and let
$$H_\alpha(x) = \log 
\left( \frac {x^{\alpha -x}\Gamma(x)}
{(x+a)^{\alpha - x  -a}\Gamma(x+a)} \right).$$
Then $-H_\alpha{'}(x)$ is completely monotonic on $(0,\infty)$ for 
$\alpha \le 1/2$, $H_\alpha{'}(x)$ is completely monotonic on
$(0,\infty)$ for $\alpha \ge 1$ and neither is completely
monotonic, for $1/2 < \alpha < 1$.  
}

The Corollary introduces ratios of gamma functions
whose asymptotic behavior \cite[6.1.47]{a&s}
\begin{equation}
z^{b-a}\frac{\Gamma(z+a)}{\Gamma(z+b)} \sim 1 + \frac{(a-
b)(a+b+1)}{2z} + \dots
\end{equation}
suggests that we look at the possible completely monotonic
character of such ratios for both gamma and $q$-gamma
functions.

Our first result is suggested by \cite[Theorem 3]{Bu:Is} on 
$(x+c)^{a-b} \Gamma(x+b)/\Gamma(x+a)$:

{\bf Theorem 2.5.} {\sl
Let $a< b \le a+1$ and let
\begin{equation}
g(x) := \left[\frac{1 - q^{x+c}}{1-q}\right]^{a-b}
\frac{\Gamma_q(x+b)}{\Gamma_q(x+a)}.
\end{equation}
Then $-(\log g(x))'$ is completely monotonic on $(-c, \infty)$
if 
$0 \le c \le (a+b-1)/2$ and $(\log g(x))'$ is  completely monotonic
on $(- a, \infty)$  if $c \ge a \ge 0$. 
}
{\sl Remarks.}  The limiting case $q
\rightarrow 1^-$ of  this theorem  was
proved by Bustoz and Ismail \cite[Theorem
3]{Bu:Is}. (Note that in the statement of \cite[Theorem 3
(ii)]{Bu:Is} ``$x > \beta$ '' should read ``$ x > \alpha$''.) 
Although, by a translation of the variable $x$, one
could assume $a = 0$ and so express the theorem in a simpler form,
we prefer to retain the $a$ and $b$ for reasons of symmetry. 
An interesting consequence of the last assertion of the theorem is
that neither $h'$ or $-h'$ is completely monotonic when 
$$c =  -\frac12 + \sqrt{ab +\frac14},$$ and hence $ (a+b-1)/2 
<c <a$.   Thus although it is true that 
\begin{equation} 
\frac{d}{dx} \log \left[ \frac{\Gamma(x+a)}{\Gamma(x+b)}\left(
x-\frac12 + \sqrt{ab +\frac14} \right)^{b-1} \right]
\end{equation}
is positive (a result motivated by the right--hand inequality in
(1.4), and proved essentially by the method given for
the case $b=1$
in \cite[Theorem 8]{Bu:Is}), it is not completely monotonic.
\newpage
{\sl Proof of Theorem 2.5}.  We note that 
\begin{eqnarray*}
\frac{d}{dx} \log g(x)& =& (b-a) \frac{q^{x+c} \log q}{1- q^{x+c}}
+ \psi_q(x+b) -\psi_q(x+a)\\
&=& - \int_0^\infty {e^{-xt}}  \left[\frac{ e^{-bt} - e^{-at}} 
{1-e^{-t}}
 + (b-a)e^{-ct} \right] d \gamma_q(t),
\end{eqnarray*}
where $\gamma_q$ is given by (1.12). 
Now we will use Lemma 1.2,  to show that the quantity in square
brackets has the appropriate sign.  In case $ c \le 
(a+b-1)/2$, we find that the integrand is 
$e^{-(x+c)t}$ times a function of $t$ which exceeds
$$ b-a - \frac{{\rm sinh}[(b-a)t/2]} {{\rm sinh}(t/2)} $$ and this
is positive by Lemma 1.2.  When  $c \ge a$, the integrand is 
$e^{-(x+a)t}$ times a function of $t$ which is less than 
$$ (b-a)e^{(b-a-1)t/2} - \frac{{\rm sinh}[(b-a)t/2]} {{\rm
sinh}(t/2)} $$ and this is negative by Lemma 1.2.
  On the other hand, if $ (a+b-1)/2 < c < a$, the quantity in
square brackets is
positive for $t$ close to $0$ and negative for large $t$. \hfill $\blacksquare$ 

The ranges of $a,b$ and $c$ in Theorem 2.5 fail to cover some
interesting cases.  For example, when $ b=c=0$, it gives
$-\log g(x)$ completely monotonic for $ a \le 0$, but gives no
information for $a >0$.  Hence it is of interest to record the
following result:

{\bf Theorem 2.6.} {\sl
Let $ 0 < q < 1$ and let 
\begin{equation}
h(x) =\log\left[ \left(
\frac{1 -q^x}{1-q} \right)^a
\frac{ \Gamma_q(x)}{\Gamma_q(x+a)} \right].
\end{equation}
Then $h'(x)$ is completely monotonic on $(0,\infty)$ for $a \ge 1$.
}

{\bf Proof:} 
From (1.11), we have
\begin{equation}
\begin{array}{ccl}
h'(x) & = & \int_0^\infty e^{-xt}\left[ \frac{e^{-at} - 1}
{1-e^{-t}} +a \right]d\gamma_q(t)\\
& = & \int_0^\infty e^{-xt + t(1-a)/2} \left[
ae^{(a-1)t/2} - \frac{{\rm sinh}(at/2)} {{\rm sinh}(t/2)} 
\right] d\gamma_q(t)
\end{array}
\end{equation} and we see from Lemma 2.2 that the quantity in
square brackets
is positive for $a >1 $.\hfill $\blacksquare$

\newpage

\section{ PSI AND $q$--PSI FUNCTIONS}

\setcounter{equation}{0}

 The psi function has particularly simple monotonicity
properties. For example, it follows from (1.8) that 
\begin{equation}
\psi'(x) = \sum_{n=0}^\infty (n+x)^{-2}, 
\end{equation}
so $\psi'(x)$ is completely monotonic on $(0,\infty)$. G.
Ronning \cite{ronning} showed that, for $0 < \alpha < 1$, $\psi'(x)
- \alpha \psi'(\alpha x) < 0,\; 0 < x < \infty$.  More generally,
we show:

{\bf Theorem 3.1.} {\sl
Let $0 < \alpha < 1,\;0 < q < 1$. Then the function
$\psi_q(x) -
\psi_{q^{1/\alpha}}(\alpha x)$ is completely monotonic,
 i.e.,
\begin{equation}
(-1)^n [\psi_q^{(n)}(x) - \alpha^n \psi_{q^{1/\alpha}}^{(n)}
(\alpha
x) ] >0, \;
0 < x < \infty, \;n=0,1,\dots. 
\end{equation}
}
{\bf Proof:} 
 Using (1.9) and (1.11) we see that 
\begin{equation} 
\psi_q(x) - \psi_{q^{1/\alpha}}(\alpha x) = 
\frac{1}{\alpha} \int_0^\infty {e^{-xt}}
\left[-\frac{\alpha}{1-e^{-t}} + \frac{1}{1-e^{-t/\alpha}} \right]
d\gamma_q(t).
\end{equation}
The quantity in square brackets is seen to be equal to
$$ \frac
{\alpha e^{(1-1/\alpha)t/2}}
{ 1 - e^{-t/\alpha}}
\left[\frac{1}{\alpha}e^{(t/2)(1/ \alpha -1)} -\frac
{{\rm sinh} (t/(2\alpha))} {{\rm sinh} (t/2)} \right],
$$
and the term in square brackets is seen to be 
 positive on using the
left--hand inequality in Lemma 1.2 (with $\alpha$ replaced by
$1/\alpha$).  The result follows.  \hfill $\blacksquare$

Gautschi's and Kershaw's inequalities (1.2) and
(1.3) suggest that we consider ratios
of the form 
\begin{equation} 
\frac{\Gamma_q(x+a)}{\Gamma_q(x+b)}\exp\left[(b-a)\psi_q(x+c)
\right].
\end{equation}
We have the following result.

{\bf Theorem 3.2.} {\sl
Let $0 < a < b,\; 0 < q \le 1$ and let $h(x)$ denote
the logarithm of the function in (3.4). Then, if $ c \ge
(a+b)/2$, 
$-h'(x)$ is completely monotonic on $(-a, \infty)$ and if $ c \le
a$, $h'(x)$ is
completely monotonic on $(-c,\infty)$.  Neither $h'$ or $-h'$ is
completely monotonic for $a < c < (a+b)/2$.
}

Bustoz and Ismail have this result in the case $q=1,\;b=1,\;
c=(a+b)/2$.

{\bf Proof:} 
 We have, from (1.9) and (1.11), 
\begin{equation}
h'(x) = - \int_0^\infty
\frac{e^{-(x+(a+b)/2)t}}{1-e^{-t}}
\left[2{\rm sinh} \frac{b-a}{2}t - (b-a)te^{((a+b)/2-c)t}
\right] d\gamma_q(t) .
\end{equation}
Using the inequality ${\rm sinh} \theta > \theta,\;\theta > 0$, we
see that the quantity in square brackets is positive when $c \ge
(a+b)/2$. When $ c \le a$, the quantity in square brackets is
\newpage
\noindent
negative on account of the inequality $ e^{-2x } > 1-2x,\;x>0$. On
the other hand, when  $a<c < (a+b)/2$ this quantity is
negative for small $t$ and positive for large $t$.
This completes the proof.
An interesting consequence of the last assertion of the theorem is
that neither $h'$ or $-h'$ is completely monotonic when $c$ is the
geometric mean of $a$ and $b$.  Thus although it is true that 
\begin{equation} 
\frac{d}{dx} \log \left[
\frac{\Gamma(x+a)}{\Gamma(x+b)}\exp\left[(b-a)\psi(x+\sqrt{
ab})
\right] \right]
\end{equation}
is positive (the result is suggested by the left--hand inequality
in (1.3) and the proof is essentially that given for the
case $b=1$
in \cite[Theorem 8]{Bu:Is}), neither it nor its negative is
completely monotonic. \hfill $\blacksquare$

A $q$--analogue of inequalities in (1.3) and
(1.4) runs as follows:

{\bf Theorem 3.3.} {\sl For $0 < q \le 1$, we have
\begin{equation}
\left(\frac{1-q^{x+s/2}}{1-q}\right)^{1-s}
< \frac{\Gamma_q(x+1)}{\Gamma_q(x+s)} < \exp[(1-s)\,
\psi_q(x+(s+1)/2)],\quad 
0 < s < 1, 
\end{equation}
where the left--hand inequality holds for $ x > -s/2$ and the
right--hand one holds for $ x > -s$.
}

{\bf Proof:} Theorem 2.5 shows that 
$$\left[\frac{1 - q^{x+s/2}}{1-q}\right]^{s-1}
\frac{\Gamma_q(x+1)}{\Gamma_q(x+s)}$$
decreases on $(-s/2,\infty)$ to its limiting value $1$. This gives
the left--hand part of (3.7).
Theorem 3.3 shows that 
$$ \exp(1-s) \psi_q\left(x +\frac{s+1 }{2} \right)
\frac{\Gamma_q(x+s)}{\Gamma_q(x+1)}$$
decreases on $(-s,\infty)$ to its limiting value $1$. This gives
the
right--hand part of (3.7).

A result of H. Alzer \cite[Theorem 1]{alzer} suggests dealing with
products of $q$--gamma functions and exponentials of {\sl
derivatives}
of $q$--psi functions.  In this connection, we prove:

{\bf Theorem 3.4} {\sl
Let $0 < q <1$, $ 0 < s < 1$ and
$$g_\alpha(x) = (1-q)^x(q^{-x} -1)
\Gamma_q(x) \exp\left[F(q^x)/\log q - \frac{1}{12} \psi_q'(x +
\alpha) \right]$$
where $F$ is given by (2.1). 
Then  $(\log g_\alpha)'$ is completely monotonic on $(0,\infty)$
for $\alpha \ge 1/2$, $-(\log g_\alpha)'$ is completely monotonic
on
$(0,\infty)$ for $\alpha \le 0$. }

\newpage

{\sl Proof.}  We have 
$$\frac{d}{dx} \log g_\alpha(x) = - \int_0^\infty e^{-xt}
p_\alpha(t) d\gamma_qt$$
where
$$p_\alpha(t) = \frac{12 - t^2 e^{-\alpha t}}{ 12(1- e^{-t})} -
\frac12 - \frac{1}{t}
$$
Now, for $\alpha \ge 1/2$, we have $p_\alpha(t) > 0$ for $t > 0$
\cite[p. 339]{alzer}.  Also when $\alpha \le 0$, we have
$p_\alpha(t) < 0$ for $t > 0$.  Thus we get the required complete monotonicity
properties of $g_\alpha$. 

If we combine this Theorem with Lemma 1.3 we get the following
extension of \cite[Theorem 1]{alzer} to which it reduces when $q
\rightarrow 1^-$.

{\bf Corollary 3.5.} {\sl
Let $0 < q <1$, $ 0 < s < 1$ and
\begin{eqnarray*}
f_\alpha(x)& =& g_\alpha(x+s)/g_\alpha(x+1)\\
&=& \left[
\frac{(1-q)^{s-1}(1-q^{x+s})^{1/2}\Gamma_q(x+s)}{(1-q^{x+1})^{1/2}
\Gamma_q(x+1)}
\right.\\
 & & \left. \times \exp
\left\{ 
\frac{F(q^{x+s}) - F(q^{x+1})}{\log q} +\frac{1}{12}
[\psi_q'(x+1+\alpha) -\psi_q'(x+s+\alpha)] 
\right\} 
\right].
\end{eqnarray*}
Then  $(\log f)'$ is completely monotonic on $(0,\infty)$ for
$\alpha \ge 1/2$, $-(\log f)'$ is completely monotonic on
$(0,\infty)$ for $\alpha \le 0$, and neither is completely
monotonic on $(0,\infty)$ for $0 < \alpha < 1/2.$
}

The limiting case $q \rightarrow 1^-$ of this last Corollary runs
as follows:

{\bf Corollary 3.6.} {\sl
Let $ 0 < s < 1$ and
\begin{eqnarray*}
f_\alpha(x) &= &
\frac{(x +1)^{x +1/2}\Gamma_q(x+s)}{
(x +s)^{x + s -1/2}
\Gamma_q(x+1)} \\ & & \times
\exp \left\{ 
s -1 +\frac{1}{12}
[\psi'(x+1+\alpha) -\psi'(x+s+\alpha)] 
\right\}.
\end{eqnarray*}
Then  $(\log f)'$ is completely monotonic on $(0,\infty)$ for
$\alpha \ge 1/2$, $-(\log f)'$ is completely monotonic on
$(0,\infty)$ for $\alpha \le 0$, and neither is completely
monotonic on $(0,\infty)$ for $0 < \alpha < 1/2.$
}

This recovers, in a slightly extended form, the main result
\cite[Theorem 1]{alzer} of H. Alzer. 
\newpage

\section{ FURTHER PRODUCTS AND QUOTIENTS}

\setcounter{equation}{0}

Results of Bustoz and Ismail \cite[Theorem 6]{Bu:Is} concerning the
ratio
$$ \frac{\Gamma(x+a)\Gamma(x+b)}
{\Gamma(x)\Gamma(x+a+b)},
$$
 suggest the  consideration of ratios
\begin{equation}
\frac{\Gamma_q(x+a_1)\Gamma_q(x+a_2)\dots\Gamma_q(x+a_n)}
{\Gamma_q(x+b_1)\Gamma_q(x+b_2)\dots\Gamma_q(x+b_n)},
\end{equation}
where
\begin{equation} a_1 + a_2 +\dots+a_n = b_1 + b_2 +\dots+b_n =s.
\end{equation}
If we treat the $a_i$ as given, each  choice of the
$b_i$ may be thought of as a partition of the sum of the $a_i$. 
An extreme case would consist of taking  the $b_i$  equal to each
other; another extreme case would be to take
all except one of the $b_i$ to be $0$. Both of these lead to
completely monotonic functions:

{\bf Theorem 4.1.} {\sl Let $a_1,\dots,a_n$ be positive numbers,
let
$n\overline{a} = a_1 +\dots a_n$, and $0 < q \le  1$. Then both
\begin{equation} 
 -\frac{d}{dx}
\log\frac{\Gamma_q(x+a_1)\Gamma_q(x+a_2)\dots\Gamma_q(x+a_n)}
{\Gamma_q(x+\overline{a})^n}
\end{equation}
and 
\begin{equation} 
\frac{d}{dx}
\log\frac{\Gamma_q(x+a_1)\Gamma_q(x+a_2)\dots\Gamma_q(x+a_n)}
{\Gamma_q({x})^{n-1}\Gamma_q(x+a_1+a_2+\dots+a_n)}
\end{equation}
are completely monotonic on $(0,\infty)$. 
}

{\bf Proof:} 
Using (1.9) and (1.11), we find that
these expressions may be written
\begin{equation}
\int_0^\infty \frac{e^{-xt}}{1-e^{-t}} \left[ ne^
{-\overline{a}t} -e^{-a_1t} - \dots -e^{-a_nt} \right]d\gamma_q(t) 
\end{equation}
and
\begin{equation}
\int_0^\infty \frac{e^{-xt}}{1-e^{-t}} \left[ n-1 +e^
{-{(a_1 + \dots a_n)}t} -e^{-a_1t} - \dots -e^{-a_nt} \right]
d\gamma_q(t) 
\end{equation}
In the first case the quantity in square brackets is positive by
Jensen's theorem while in the second case its positivity follows
from 
\begin{equation} n-1 +z_1z_2\dots z_n - (z_1 + z_2 + \dots +z_n)
\ge 0, \;\;0
\le z_i
< 1.
\end{equation}
This is clear since it may be established by induction on $n$ that 
\begin{equation}
n-1 +z_1z_2\dots z_n - (z_1 + z_2 + \dots +z_n)
= \sum_{j=2}^n(1-z_j)(1-z_1z_2\dots z_{j-1}).
\end{equation} \hfill $\blacksquare$
\newpage
The positivity of the quantity (34) leads to the
inequality
$$\Gamma_q(1+a_1)\Gamma_q(1+a_2)\dots\Gamma_q(1+a_n) \ge
\Gamma_q(1+\overline{a})^n,\; a_i > 0, 
$$
which is known already in the case $q=1$ \cite[p. 285]{Mi}.

\section{BOUNDS IN THE COMPLEX PLANE}
\setcounter{equation}{0}

A generalization of the Phragm\'en--Lindel\"of theorem is used in 
\cite[pp. 68--70]{rademacher} to show that 
\begin{equation}
\left|\frac{\Gamma(s+c)}{\Gamma(s)} \right| \le |s|^c, \;0 \le c
\le 1, \; {\rm Re}(s)  \ge (1-c)/2.
\end{equation}
Here we show that the same method can be used to get bounds for
more complicated functions involving gamma functions.

The Phragm\'en--Lindel\"of theorem runs as
follows \cite[p. 59]{rademacher}:

{\bf Theorem 5.1.} {\sl
Let $f(z)$ be analytic in the strip $S(\alpha, \beta) =
\{z| z= x+iy,\; \alpha < x < \beta \}$. Let us assume $|f(z)| \le
1$ on the boundaries $x = \alpha$
and $x = \beta$ and moreover,
$$ |f(z)| < Ce^{e^{k|y|}} $$
for some $C > 0$ and $ 0 < k < \pi/(\beta - \alpha)$.
Then $|f(z)| \le 1$ throughout the strip $S(\alpha,\beta)$.
}

By Theorem 4.1, the function
\begin{equation}
f(x) := \frac{\Gamma(x+a)\Gamma(x+b)}
{\Gamma(x)\Gamma(x+a+b)},\quad a, b \ge 0.
\end{equation}
is increasing on the interval $(0, \infty)$  to its limit $1$.
 Hence we have 
\begin{equation}
\left| \frac{\Gamma(x+a)\Gamma(x+b)}
{\Gamma(x)\Gamma(x+a+b)} \right| \le 1,\; x > 0.
\end{equation}
We have, in fact:

{\bf Theorem 5.2.} {\sl
We have, for $0 \le a \le 1$, $b \ge 0$,
\begin{equation}
\left| \frac{\Gamma(s+a)\Gamma(s+b)}
{\Gamma(s)\Gamma(s+a+b)} \right| \le 1, \; {\rm Re} \;s > 
\frac{1-a-b}{2}.
\end{equation}
}

{\sl Proof.} The method follows the proof of \cite[Theorem
A, p. 68]{rademacher}.
Since the assertion is trivial for $a=0$, we may as well choose $a
>0$.  We choose a complex number $s= \sigma + i \tau$ satisfying
the hypotheses of the theorem and let 
$$ f(z) = \frac{\Gamma(a+2\sigma -z)\Gamma(b+z)}
{\Gamma(z)\Gamma(a+b + 2 \sigma -z)}.
$$
\newpage
Clearly 
$$ |f(s)| = \left| \frac{\Gamma(s+a)\Gamma(s+b)}
{\Gamma(s)\Gamma(s+a+b)} \right| $$
so we have to show that $|f(s)| \le 1$. 
Now with $\alpha = (a-1)/2 +\sigma,\;\;\beta = a/2 +\sigma$, and
using $\overline{\Gamma(z)} = \Gamma(\overline{z})$, we get
$$|f(\alpha+it)| = \frac{| \alpha +it|}
{|b+ \alpha +it|} \le 1$$ where we have used 
${\rm Re}\; s \ge (1-a-b)/2$.
Also $$|f(\beta+it)| = 1,$$
and (20) shows that the growth condition of Theorem 5.1 is
satisfied.  Using this theorem, we find that
$|f(z)| \le 1$, for $\alpha \le {\rm Re}\;z \le\beta$, so, in
particular, $|f(s)| \le 1$. 

{\bf Acknowledgments.} The authors' work was supported partially by
NSF grant DMS 9203659 and by NSERC Canada grant A5199.  We are
grateful to Ruiming Zhang for pointing out the relevance of
reference \cite{rademacher}.  We thank a referee for corrections
and suggestions.

\end{document}